\documentclass[smallcondensed]{svjour3}
\usepackage[numbers]{natbib}

\usepackage[utf8]{inputenc} % For accents é, ö, etc.
\usepackage{todonotes}
\usepackage{subcaption}
\usepackage{stmaryrd} % Brackets for v_jump
\usepackage{amssymb}
\usepackage[makeroom]{cancel} % To cancel terms in eqs
\usepackage{bookmark} % To add bookmarks in pdf
\usepackage{empheq}
\usepackage{ulem}
\usepackage{tcolorbox}  % For todo boxes
\usepackage{pdflscape} % For landscape pages
\usepackage{tikz-cd}

\usepackage{qting}
\newcommand{\change}[1]{{\leavevmode\color{black}#1}}

\newcommand{\bigpartialderiv}[2]{ \frac{\partial {#1}}{\partial {#2} } }

\begin{document}

\let\WriteBookmarks\relax
\def\floatpagepagefraction{1}
\def\textpagefraction{.001}

\title{Mimetic Metrics for the DGSEM}
\author{Daniel Bach \and Andrés Rueda-Ramírez \and David A. Kopriva \and Gregor J. Gassner}
\institute {Daniel Bach (University of Cologne, \email{daniel.bach@uni-koeln.de})   \and Andrés Rueda-Ramírez (Universidad Politécnica de Madrid, \email{am.rueda@upm.es}) \and David A. Kopriva (Florida State University and San Diego State University, \email{dkopriva@fsu.edu}) \and Gregor J. Gassner (University of Cologne, \email{ggassner@uni-koeln.de})}

% \sep 

\maketitle
\subclass{65Mxx}

\begin{abstract}
Free-stream preservation is an essential property for numerical solvers on curvilinear grids. 
Key to this property is that the metric terms of the curvilinear mapping 
satisfy discrete metric identities, i.e., 
have zero divergence. Divergence-free metric terms are furthermore essential for entropy stability on curvilinear grids. 
We present a new way to compute the metric terms for discontinuous Galerkin spectral element methods (DGSEMs) that guarantees they are divergence-free. 
The proposed mimetic approach uses projections that fit within the de Rham Cohomology.
\end{abstract}

\keywords{Free-Stream Preservation; Mimetic Methods; Discontinuous Galerkin; Divergence Free Methods; Curved Meshes}

\section{Introduction}

\q{c16_r2}{
We focus on discontinuous Galerkin spectral element methods (DGSEMs) for \change{systems of conservation laws} of the form \change{
\begin{equation}
    \label{eq:con_law}
    \frac{\partial}{\partial t} u_i + \text{div}_x\left(\textbf{F}_i(\mathbf{u})\right) = 0 \qquad i = 1,\dots,m
\end{equation}}
on a domain $\Omega \times [0,T]$, \change{where $\mathbf{u} = \left(u_1, \dots, u_m\right)^T$ is the vector of conserved variables, $\textbf{F}_i(\textbf{u})$ is the vector-valued flux function of the $i$-th equation and $\text{div}_x$ denotes the divergence operator in physical space. 
In this paper, we denote differential operators in physical space with a subscript $x$, as in ``div$_x$''. 
Furthermore, we use a bold font for vector-valued variables and functions, while we do not use bold letters for scalar quantities.
In the following, we will only look at one of the equations at a time, so we will omit the subscript, $i$. }}

\q{c16_r2_2}{
To solve the \change{system} numerically we partition the domain $\Omega$ with a conforming mesh composed of quadrilateral (in 2D) or hexahedral (in 3D) elements.
\change{Let $\left(x, y, z\right)^T = \left(x_1, x_2, x_3\right)^T = \textbf{x}: E \longrightarrow Q$ be the mapping from the reference element $E$ to an element $Q$ in the physical domain. The computational coordinates are denoted as $\boldsymbol{\xi} = \left(\xi^1,\xi^2,\xi^3\right)^T = \left(\xi,\eta,\zeta\right)^T \in [-1, 1]^3$.
For clarity, we will omit the subscript $x$ when differential operators are applied with respect to the reference element coordinates.}
}
After \change{\eqref{eq:con_law}} is transformed to the reference element, it becomes
\begin{equation}\label{eq:pde_tranfo}
     \frac{\partial}{\partial t} \left(Ju\right) + \text{div}\left((J \textbf{a}^s \cdot \textbf{F}(\mathbf{u}))_{s=1}^3\right) = 0,
\end{equation}
where $J$ is the mapping Jacobian and $J\textbf{a}^i = J \textbf{grad}_x (\xi^i)$ are the metric terms,
\begin{equation}\label{eq:crossproduct}
    \begin{split}
    J\textbf{a}^i &= \bigpartialderiv{\textbf{x}}{\xi^j}
    \times
    \bigpartialderiv{\textbf{x}}{\xi^k} ,\,\,
    (i,j,k)\text{ cyclic}, \\
    \text{ or equivalently \cite{kopriva}}, J\textbf{a}^i_n &= - \hat{\textbf{x}}_i \cdot \textbf{curl} \left( x_m \textbf{grad} \left(x_l\right)\right),\,\,
    (n,m,l)\text{ cyclic},
    \end{split}
\end{equation}
Here $\hat{\textbf{x}}_i$ are the canonical space unit vectors.

If we consider a constant state $\mathbf{u}=\mathbf{c}$ and its corresponding %constant 
flux $\textbf{C} = \textbf{F}(\mathbf{c}) = const$, with constant or periodic boundary conditions, we get from \eqref{eq:pde_tranfo} %the time derivative
\begin{equation}
\frac{\partial}{\partial t} \left(Ju\right) = - \text{div}\left(\left(J\textbf{a}^s \cdot \textbf{C}\right)_{s = 1}^3\right)  = 0,
\end{equation}
%where we used the fact 
using that for the analytical global mapping with the conditions specified the divergence of the metric terms is zero - the constant state stays constant in time.

For the DGSEM (see, e.g., \cite{kopriva}), we use a tensor product ansatz and interpolate $Ju$ and the flux functions $J \textbf{a}^s \cdot \textbf{F}(\mathbf{u})$ with polynomials of degree $N$ in each direction. 
Let $I^N$ be the tensor product interpolation operator and $l_{\textbf{k}}$ the tensor product Lagrange polynomials. 
\change{We multiply \eqref{eq:pde_tranfo} with test functions, %which in this case are 
using the same Lagrange polynomials, and integrate over the reference element.} Then we use integration-by-parts, replace the boundary flux with the numerical flux, use the summation-by-parts property and collocate interpolation and quadrature using Legendre-Gauss-Lobatto (LGL) nodes to arrive at the strong form of the DGSEM,
\begin{equation}
    \begin{split}
    &\frac{\partial}{\partial t} \left( Ju \right)_{\textbf{i}} + \text{div}\left( \left( I^N(J \textbf{a}^s \cdot \textbf{F}(\textbf{u})) \right)_{s = 1}^3 \right)_{ \textbf{i}}\\
    = &\frac{1}{\omega_{\textbf{i}}} \int_{\partial E,N}  \left( I^N\left(F^{\text{num}}(\mathbf{u}_l, \mathbf{u}_r, \textbf{n})\right) - \left( I^N\left(\textbf{F} \cdot \textbf{n}\right) \right) \right) l_{\textbf{i}} \text{ d}\xi,
    \end{split}
\end{equation}
for all multi-indices $\textbf{i} \in \{0, \dots, N\}^n$. 
Here the multiindex $\textbf{i}$ denotes the value of the function at the corresponding grid point, $E$ is the reference element, $\omega_{\textbf{i}}$ is the product of the LGL quadrature weights at the corresponding grid point and $\int_{\partial E, N}$ denotes the numerical LGL-quadrature of the boundary integral (see \cite{kopriva}). 
\q{c17_r2}{\change{Finally, $\textbf{n}$ denotes the outward directed normal vector on the interface and $u_l$ and $u_r$ denote the left- and right-sided value on the interface, where left always means the value from `inside' the element under consideration.}}

Substituting the constant state $\mathbf{u}=\mathbf{c}$ into the DGSEM, the boundary term cancels due to the consistency of the numerical flux leaving
\begin{equation}
    \frac{\partial}{\partial t} \left( Ju \right)_{\textbf{i}} = -  \text{div} \left(\left(I^N\left(J \textbf{a}^s \cdot \textbf{C}\right) \right)_{s=1}^3\right)_{ \textbf{i}},
\end{equation}
which is zero if the divergence of the discrete metric terms is zero,
\begin{equation} \label{eq:discretemetricidentities}
    \text{div} \left(\left(I^N\left(J\textbf{a}_i^s \right) \right)_{s = 1}^3 \right) = 0 \hspace{1cm} \forall i \in \{1,\dots,3\}.
\end{equation}

Even if the mapping $\textbf{x}$ is polynomial, the metric terms $J\textbf{a}^i_n$ computed using \eqref{eq:crossproduct} are of a higher polynomial order.
For the operator $I^N$ to still act as the identity operator in \eqref{eq:discretemetricidentities}, we would have to restrict the mapping order to $q \le \frac{N}{2}$. 
A different approach is used in \cite{kopriva} where the metric terms are computed as
\begin{equation}
    \label{eq:3dmetricdiscrete}
    J\textbf{a}^i_n \approx -\frac{1}{2}\hat{\textbf{x}}_i \cdot \textbf{curl}\left(\mathbf{I}^N\left( {x}_l \textbf{grad}\left({x}_m\right) - {x}_m \textbf{grad}\left({x}_l\right)\right) \right)
    \qquad
    (n,m,l) \text{ cyclic}.
\end{equation}
Here $\mathbf{I}^N$ is the component-wise interpolation.
Since \eqref{eq:3dmetricdiscrete} takes the analytical curl of the interpolation and the divergence of a curl is zero, the curl form approximation of the metric terms guarantees free-stream preservation.

\section{Mimetic Approach}

\q{c1_r1}{
\change{We propose an alternative to} the curl form of \cite{kopriva} using mimetic projections coming from finite element exterior calculus \cite{feec}.}
\q{c6_r2}{
\change{To apply this framework, we look at the de Rham cochain complex in 3D, and need to choose for each space a finite dimensional subspace $V_i$ and a projection $p^i$ so that the diagram \ref{fig:derham3d} commutes.
\begin{figure}[h]
    \centering
    \begin{tikzpicture}[node distance = 1.0cm]
        \label{derham3d}
    	\node (1) at (2,0) {$H^1$};
    	\node (2) at (4,0) {$H\left(\textbf{curl}\right)$};
    	\node (3) at (6,0) {$H\left(\text{div}\right)$};
    	\node (4) at (8,0) {$L^2$};
    	\node (7) at (2,-1.5) {$V_1$};
    	\node (8) at (4,-1.5) {$V_2$};
    	\node (9) at (6,-1.5) {$V_3$};
    	\node (10) at (8,-1.5) {$V_4$};
    	\draw[->] (1) -- node [midway, above] {\textbf{grad}} (2);
    	\draw[->, red] (2) -- node [midway, above] {\textbf{curl}} (3);
    	\draw[->] (3) -- node [midway, above] {div} (4);
    	\draw[->] (7) -- node [midway, above] {\textbf{grad}} (8);
    	\draw[->, blue] (8) -- node [midway, above] {\textbf{curl}} (9);
    	\draw[->] (9) -- node [midway, above] {div} (10);
      	\draw[->] (1) -- node [midway, right] {$p^1 = I^N$} (7);
    	\draw[->, blue] (2) -- node [midway, right] {$\textbf{p}^2$} (8);
    	\draw[->, red] (3) -- node [midway, right] {$\textbf{p}^3$} (9);
    	\draw[->] (4) -- node [midway, right] {$p^4$} (10);
    \end{tikzpicture}
    \caption{Diagram of the continuous and discrete 3D de Rham complex}
    \label{fig:derham3d}
\end{figure}
}}
\change{
To choose spaces and projections that fit our needs, we follow the work of Gerritsma et al. \cite{gerritsma,gerritsma_mimetic_framework} starting from the 1D de Rham complex \ref{fig:derham1d}.
\begin{figure}[h]
    \centering
    \begin{tikzpicture}[node distance = 1.0cm]
    	\node (1) at (2,0) {$H^1$};
    	\node (2) at (4,0) {$L^2$};
    	\node (7) at (2,-1.5) {$V_1$};
    	\node (8) at (4,-1.5) {$V_2$};
    	\draw[->] (1) -- node [midway, above] {$\bigpartialderiv{}{\xi}$} (2);
    	\draw[->] (7) -- node [midway, above] {$\bigpartialderiv{}{\xi}$} (8);
      	\draw[->] (1) -- node [midway, right] {$p^1$} (7);
    	\draw[->] (2) -- node [midway, right] {$p^2$} (8);
    \end{tikzpicture}
    \caption{Diagram of the continuous and discrete 1D de Rham complex}
    \label{fig:derham1d}
\end{figure}
}
\change{
We choose for each element the spaces $V_1 = \mathbb{P}^{N}$ and $V_2 = \mathbb{P}^{N-1}$, which are the polynomial spaces of degree $N$ and $N-1$ respectively. 
To construct the projections we impose a sub-grid of the LGL nodes in each element, which gives an irregular distribution of sub-intervals between those nodes. 
Then we define basis functions on $V_1$ and $V_2$, where the $V_1$ basis corresponds to the LGL nodes and the $V_2$ basis corresponds to the sub-intervals. 
For $V_1$ these basis functions are just the Lagrange polynomials $(l_i)_{i = 0}^N$, while for $V_2$ they are the so-called edge polynomials, e.g. \cite{gerritsma}, $(h_i)_{i = 1}^N$ defined as}
\begin{equation}
    \label{eq:deriv}
    h_i(\xi) = -\sum_{j = 0}^{i-1} \bigpartialderiv{l_j}{\xi}(\xi),
\end{equation}
\change{with the property
\begin{equation}
    \label{eq:hist_property}
    \int_{\xi_{j-1}}^{\xi_j} h_i(\xi) \text{d}\xi = \delta_{i,j},
\end{equation}
where $\delta_{i,j}$ is the Kronecker delta. Eq. \eqref{eq:hist_property} is a property analogous to that of the Lagrange polynomials, but for the sub-interval integrals instead of the point values.} 
\q{c12_r2}{\change{Using these basis functions we can define in each element the projections $p^1$ and $p^2$ by setting $p^1 = I^N$
and $p^2$ to be the histopolation, %in each element
which is the approximation by polynomials preserving the sub-interval integrals. 
For one element this means
\begin{equation}
    p^2(f) = \sum_{i = 1}^N \left(\int_{\xi_{i-1}}^{\xi_i} f(\xi) \text{d}\xi\right) h_i.
\end{equation}
}}
\change{As shown in \cite{gerritsma}, the 1D diagram  \ref{fig:derham1d} commutes with these definitions. Since the LGL nodes include the boundary points of the element, we get inter-element continuity for $V_1$ while there is no such requirement for $V_2$.

For the 3D case we can now use a tensor product ansatz of Lagrange and edge basis functions, and define the spaces as 
\begin{align}
V_1 &= \mathbb{P}^{N} \otimes \mathbb{P}^{N} \otimes \mathbb{P}^{N},\\
V_2 &= \left(\mathbb{P}^{N-1} \otimes \mathbb{P}^{N} \otimes \mathbb{P}^{N}, \mathbb{P}^{N} \otimes \mathbb{P}^{N-1} \otimes \mathbb{P}^{N}, \mathbb{P}^{N} \otimes \mathbb{P}^{N} \otimes \mathbb{P}^{N-1}\right)^T,\\
V_3 &= \left( \mathbb{P}^{N} \otimes \mathbb{P}^{N-1} \otimes \mathbb{P}^{N-1}, \mathbb{P}^{N-1} \otimes \mathbb{P}^{N} \otimes \mathbb{P}^{N-1}, \mathbb{P}^{N-1} \otimes \mathbb{P}^{N-1} \otimes \mathbb{P}^{N}\right)^T,\\
V_4 &= \mathbb{P}^{N-1} \otimes \mathbb{P}^{N-1} \otimes \mathbb{P}^{N-1},
\end{align}
with continuity requirements derived analogously as in the 1D case, i.e. where $\mathbb{P}^{N}$ indicates directions with inter-element continuity, and $\mathbb{P}^{N-1}$ indicates directions that are generally discontinuous.} 
We note that these continuity requirements couple the nodal degrees of freedom similar to a continuous finite element approach.

\change{The projection operators are similarly constructed based on the tensor product approach.} 
We give the first component of $\textbf{p}^2: H(\textbf{curl})\rightarrow V_2$ as an example,
\begin{equation}
    \textbf{p}^2_1\left(\textbf{f}\right) = \sum_{i=1}^N \sum_{j = 0}^N \sum_{k = 0}^N \left(\int_{\xi_{i-1}}^{\xi_{i}} f_1(s,\eta_j,\zeta_k) \text{ ds}\right) h_i(\xi)l_j(\eta)l_k(\zeta).
\end{equation}
\change{This construction is one of many possible approaches, for example the one used by Nédélec for cubic elements \cite{nedelec}. 
The Nédélec spaces are the same, but the Nédélec elements use higher moment degrees of freedom instead of a sub-grid nodes, leading to different projections (see the discussion by Gerritsma et al. in \cite{gerritsma_mimetic_framework}}).

\change{We can now
%use the definition of the discrete spaces and projections 
define our approximative metric terms. 
%Our choices give a 
Since the 3D diagram \ref{fig:derham3d} commutes, we} have two options with which to construct the discrete metric terms that result in the same approximation (up to machine precision). 

{\color{red}Option red}: Project $\textbf{curl}\left(x_m \textbf{grad}\left(x_l\right)\right)$ directly via $\textbf{p}^3$ to get
\begin{equation}
    \label{eq:op_red}
   \textbf{p}^3 \left(\left(J\textbf{a}^i_n\right)_{i = 1}^3\right) = \textbf{p}^3(\textbf{curl}\left(x_m \textbf{grad} \left(x_l\right)\right))\in V_3.
\end{equation}

{\color{blue} Option blue}: Project $x_m \textbf{grad}\left(x_l\right)$ via $\textbf{p}^2$ to get $\textbf{p}^2\left(x_m \textbf{grad} \left(x_l\right)\right) \in V_2$. 
We compute the metric terms by applying the curl, i.e. 
\begin{equation}
    \label{eq:op_blue}
    \textbf{p}^3 \left(\left(J\textbf{a}^i_n\right)_{i = 1}^3\right) = \textbf{curl}\left(\textbf{p}^2\left(x_m \textbf{grad}\left(x_l\right)\right)\right) \in V_3.
\end{equation} %It then holds that 
\change{We then get free-stream preservation by construction:}
\begin{equation}
        \text{div}\left(\textbf{p}^3\left(\textbf{curl} \left(x_m \textbf{grad}\left(x_l\right)\right)\right)\right) = \text{div}\left(\textbf{curl}\left(\textbf{p}^2\left(x_m \textbf{grad}\left(x_l\right)\right)\right)\right) = 0
\end{equation}
%by construction, which gives free-stream preservation. 
\q{c3_r1}{
\begin{remark}
\change{
We note that for a fixed mesh, computing the metric terms is a pre-processing step with typically negligible impact on the overall compute time for 2D and 3D simulations. 
To make life simple, we could for instance first approximate the geometry with a piece-wise $C^0$ polynomial ansatz by interpolation and use the resulting approximative mapping to compute the metric terms. 
By using standard numerical integration and differentiation procedures, the computational complexity is comparable to the approach presented by Kopriva \cite{kopriva}.

It is also possible to directly use the (analytical) mapping as is. 
%Theoretically, one could directly integrate and differentiate exactly, if possible. 
If possible, one could integrate and differentiate exactly. 
As an alternative to exact evaluation of metric terms, one could use again a numerical integration (and differentiation) approach that is implemented so that the errors are of machine precision magnitude, e.g., by using sufficient quadrature points and automatic differentiation. 
%In theory, it is also possible to directly use a NURBS representation of the geometry from a CAD software, although this would introduce additional software development complexity, such as e.g. a custom CAD geometry interface.
}
\end{remark}
}

\begin{remark}
In 2D, the approach in \cite{kopriva} is to interpolate the mapping components and then take the analytic 2D curl to achieve free-stream preservation. 
This approach is equivalent to the mimetic approach because the metric terms in 2D depend linearly on the mapping components, i.e., \eqref{eq:crossproduct} reduces to 
\begin{equation}
    J\textbf{a}^1 = \left(\bigpartialderiv{}{\eta}y, -\bigpartialderiv{}{\eta}x \right)^T,\qquad J\textbf{a}^2 = \left(-\bigpartialderiv{}{\xi}y, \bigpartialderiv{}{\xi}x \right)^T.
\end{equation}
In other words, for 2D we have the smaller commuting diagram \ref{derham2d} where $\textbf{curl}_v(f) := \left(\bigpartialderiv{f}{\eta}, -\bigpartialderiv{f}{\xi}\right)^T$. 
The approach for 2D in \cite{kopriva} is equivalent to taking the blue path.
\q{c6_r2_2}{\change{
\begin{figure}[h]
    \centering
    \begin{tikzpicture}
        \label{derham2d}
		\node (1) at (2,0) {$H^1$};
		\node (2) at (4,0) {$H\left(\text{div}\right)$};
		\node (3) at (6,0) {$L^2$};
		\node (6) at (2,-1.5) {$V_1$};
		\node (7) at (4,-1.5) {$V_3$};
		\node (8) at (6,-1.5) {$V_4$};
		\draw[->, red] (1) -- node [midway, above] {$\textbf{curl}_v$} (2);
		\draw[->] (2) -- node [midway, above] {div} (3);
		\draw[->] (7) -- node [midway, above] {div} (8);
		\draw[->, blue] (6) -- node [midway, above] {$\textbf{curl}_v$} (7);
  		\draw[->, blue] (1) -- node [midway, right] {$p^1 = I^N$} (6);
		\draw[->, red] (2) -- node [midway, right] {$\textbf{p}^3$} (7);
		\draw[->] (3) -- node [midway, right] {$p^4$} (8);
	\end{tikzpicture}\\
    \caption{Diagram of the continuous and discrete 2D de Rham complex}
    \label{derham2d}
\end{figure}
}}
\end{remark}

\section{Numerical Example}

\q{c4_r1}{\change{As an example, we discretize the 3D compressible Euler equations}
\begin{equation}
    \begin{pmatrix}
        \rho \\
        \rho v_1 \\
        \rho v_2 \\
        \rho v_3 \\
        \rho e
    \end{pmatrix}_t
    + \begin{pmatrix}
        \rho v_1 \\
        \rho v_1^2 + p \\
        \rho v_1 v_2 \\
        \rho v_1 v_3 \\
        \left(\rho e + p\right) v_1
    \end{pmatrix}_x
    + \begin{pmatrix}
        \rho v_2 \\
        \rho v_1 v_2 \\
        \rho v_2^2 + p \\
        \rho v_2 v_3 \\
        \left(\rho e + p\right) v_2
    \end{pmatrix}_y
    + \begin{pmatrix}
        \rho v_3 \\
        \rho v_1 v_3 \\
        \rho v_2 v_3 \\
        \rho v_3^2 + p\\
        \left(\rho e + p\right) v_3
    \end{pmatrix}_z
    = \begin{pmatrix}
        0 \\
        0 \\
        0 \\
        0 \\
        0
    \end{pmatrix},
\end{equation}
\change{where $e$ is the specific total energy and} \begin{equation}
    p = \left(\gamma - 1\right) \left( \rho e - \frac{1}{2} \rho \left(v_1^2 + v_2^2 + v_3^2 \right) \right),
\end{equation}
\change{with $\gamma = 1.4$ and initial condition} \begin{equation*}
    \rho \equiv 1.0\text{, }\rho v_1 = 0.1 \text{, } \rho v_2 \equiv -0.2 \text{, } \rho v_3 \equiv 0.7 \text{, } \rho e \equiv 10.0.
\end{equation*}
}
We use the mapping 
%\begin{align}\label{example:mapping}
%    x &= \xi + \theta (\xi,\eta,\zeta),
%    \\
%    y &= \eta + \theta (\xi,\eta,\zeta),
%    \\
%    z &= \zeta + \theta (\xi,\eta,\zeta),
%\end{align}
\begin{equation}
    \textbf{x}(\boldsymbol{\xi}) = \boldsymbol{\xi} + \theta (\boldsymbol{\xi}) (1,1,1)^T
\end{equation}
from the reference domain $[-1,1]^3$ with
%\begin{equation}
%    \theta (\xi,\eta,\zeta) := 0.1 \cos(\pi \xi) \cos(\pi %\eta) \cos(\pi \zeta),
%\end{equation}
\begin{equation}
    \theta (\boldsymbol{\xi}) := 0.1 \cos(\pi \xi) \cos(\pi \eta) \cos(\pi \zeta),
\end{equation}
and end time $T = 1$. We use a conforming mesh with two elements in each direction and periodic boundary conditions. \q{c14_r2}{\change{For time integration we use a fourth-order, five-stage low storage Runge Kutta method by Carpenter and Kennedy (\cite{carpenter}, page 13, solution 3) with a CFL of $0.2$, which is stable for all polynomial degrees investigated.
The time step $\Delta t$ is calculated by 
\begin{equation}
    \Delta t = \text{CFL} \frac{2}{|\lambda_{\text{max}} (N+1)|},
\end{equation}
where $\lambda_{\text{max}}$ is the maximum wave speed of the PDE system evaluated at the LGL nodes scaled by the geometry.}}
Simulations are done with a slightly modified version of the DGSEM code of the Julia package Trixi.jl (see \cite{ranocha2022adaptive} and \cite{schlottkelakemper2020trixi}) performed on a Macbook Pro M2, single thread, with MacOS 14.5.

We examine both free-stream preservation and the approximation of the analytic metric terms by varying the polynomial degree up to $N = 25$ for both the curl form by Kopriva in \cite{kopriva} and the mimetic variant. Both errors are measured in the $L_2$ and the $L_{\infty}$ norm. The norms are computed numerically using LGL integration with 51 points in each direction in every element for the $L_2$ norm, and by the \change{maximum} of the point value errors on the same points as the discrete evaluation of the $L_{\infty}$ norm.
\q{c15_r2}{\change{For the metric terms we compute on each of the 51 LGL points in each direction per element the norms $\left|\left|\text{vec}\left(Ja_n^i\right)_{i,n = 1}^3\right|\right|_2$ and $\left|\left|\text{vec}\left(Ja_n^i\right)_{i,n = 1}^3\right|\right|_{\infty}$ and then apply the LGL integration formula for the $L_2$ or choose the maximum value for the $L_{\infty}$ norm respectively.}} 

\q{c4_r1_2}{
\begin{figure}
    \centering
    \begin{subfigure}[b]{0.49\linewidth}
\includegraphics[trim=0 0 0 0, clip,
width=\linewidth]{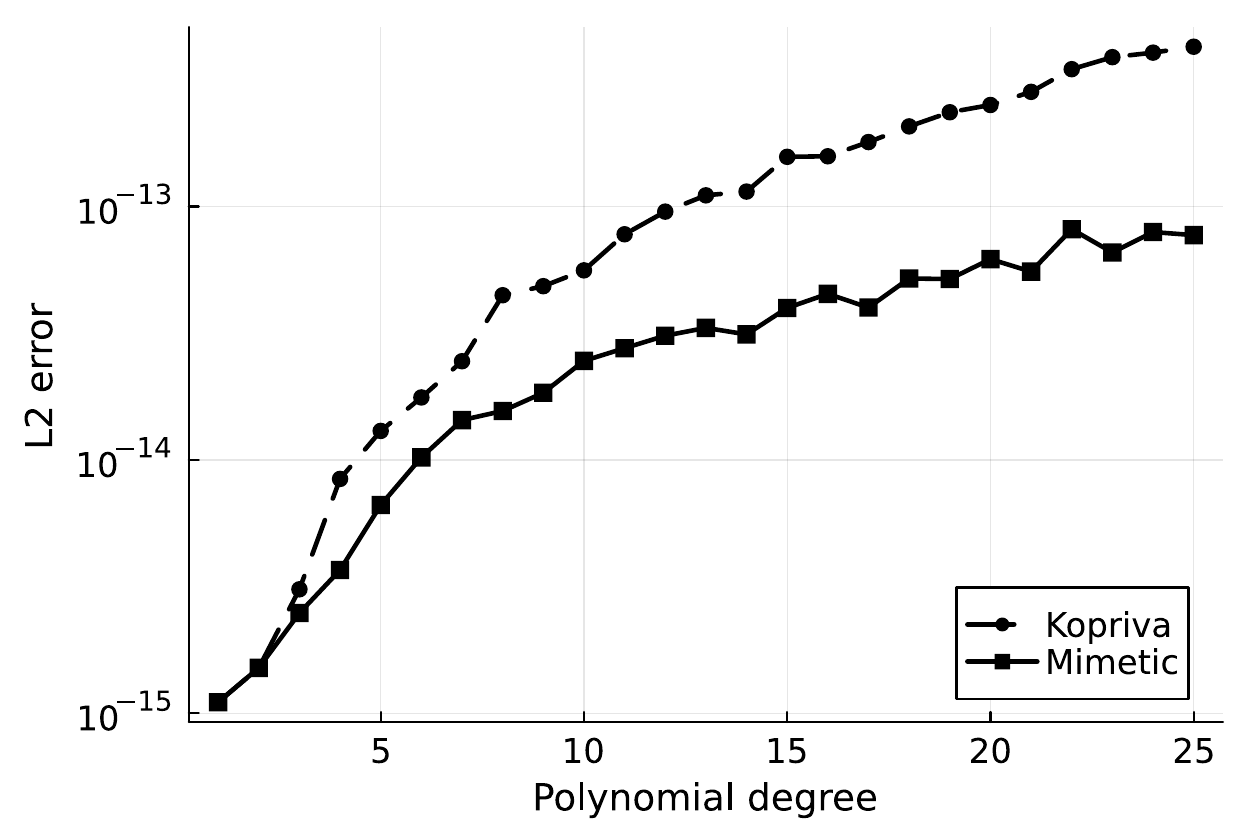}
\caption{Discrete $L_2$-error}
\label{fig:fsp_l2_error}
\end{subfigure}
\begin{subfigure}[b]{0.49\linewidth}
\includegraphics[trim=0 0 0 0, clip,
width=\linewidth]{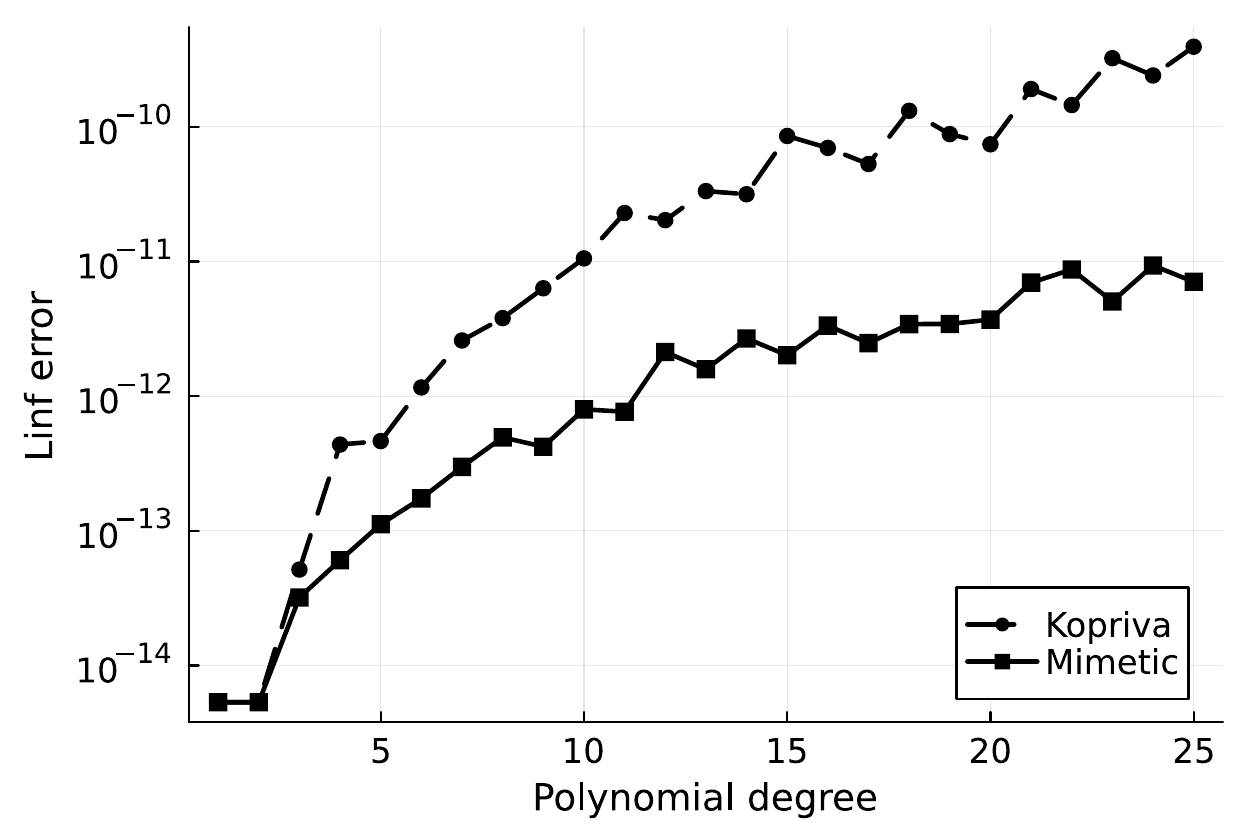}
\caption{Discrete $L_{\infty}$-error}
\label{fig:fsp_linf_error}
\end{subfigure}
\caption{Free-stream preservation errors of $\rho e$ with method of Kopriva \cite{kopriva}, Eq. \eqref{eq:3dmetricdiscrete}, and the Mimetic approach, Eq. \eqref{eq:op_blue}.}
\end{figure}
}
\begin{figure}
    \centering
    \begin{subfigure}[b]{0.49\linewidth}
\includegraphics[trim=0 0 0 0, clip,
width=\linewidth]{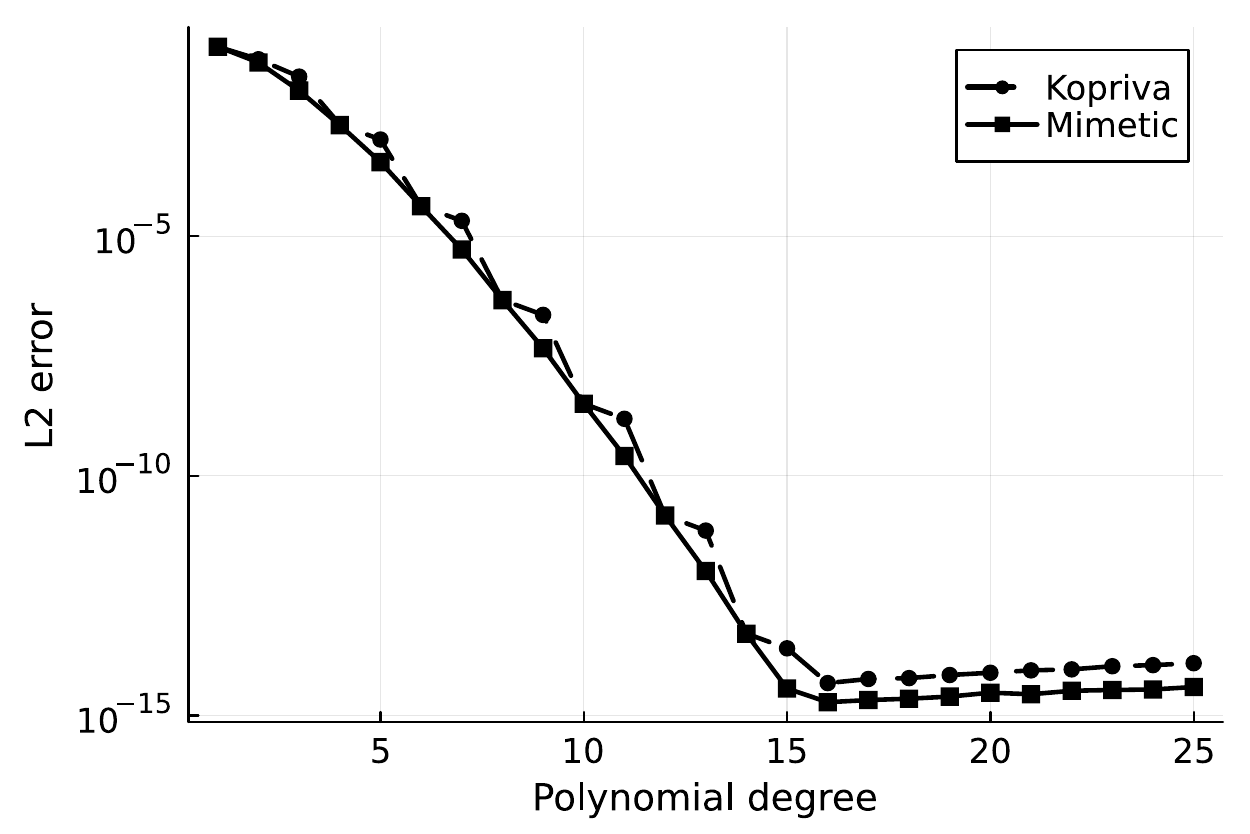}
\caption{Discrete $L_2$-error}
\label{fig:contravar_l2_error}
\end{subfigure}
\begin{subfigure}[b]{0.49\linewidth}
\includegraphics[trim=0 0 0 0, clip,
width=\linewidth]{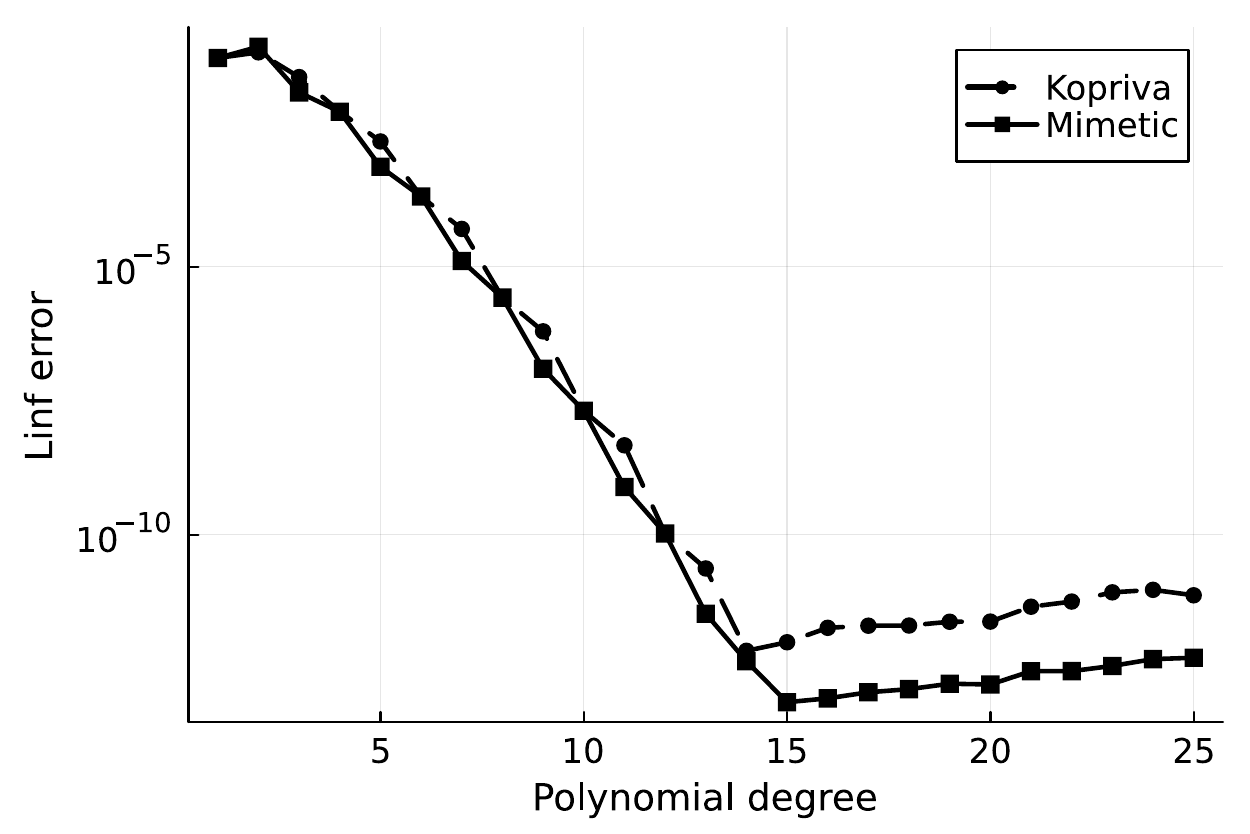}
\caption{Discrete $L_{\infty}$-error}
\label{fig:contravar_linf_error}
\end{subfigure}
\caption{Error in the \change{metric terms} with method of Kopriva, \eqref{eq:3dmetricdiscrete}, and the Mimetic approach, see \eqref{eq:op_blue}.}
\end{figure}

\q{c18_r2}{In Fig. \ref{fig:fsp_l2_error} and \ref{fig:fsp_linf_error} we can see that free-stream preservation is achieved with both the curl and the mimetic approach, and that the \change{rounding} errors grow with increasing polynomial degree for both. We also observe that the mimetic approach has \change{rounding} errors smaller than Kopriva's curl form by over an order of magnitude in the $L_{\infty}$ norm.} 
\q{c_machine_precision_r2}{In Figure \ref{fig:contravar_l2_error} and \ref{fig:contravar_linf_error} we see convergence of the \change{metric terms} to machine precision, though the \change{rounding} errors are again higher for Kopriva's approximation.
}
\section{Conclusion}
We have constructed an alternative way to compute the metric terms for the DGSEM scheme in a divergence-free way via a mimetic approach using \change{de Rham complexes} for both two and three space dimensions. 
In two space dimensions the method of \cite{kopriva} and the mimetic approach are analytically the same, but they differ in three dimensions. Both approaches do obtain free-stream preservation for conforming meshes and exhibit a similar convergence behavior for the \change{metric terms}. 
However, in all cases the mimetic approach has better rounding error properties. We remark that the finite element exterior calculus projections do not only exist for quadrilateral or hexahedral meshes but also, for example, for triangular or tetrahedral meshes \cite{feec}.

\section*{Declarations}

\textbf{Conflicts of Interest}: The authors declare that there are no conflicts of interest.

\paragraph{Data Availability}: A reproducibility repository can be accessed under \url{https://github.com/amrueda/paper_2024_mimetic_metrics}.

\paragraph{Funding}
This work was supported by a grant from the Simons Foundation (\#426393, \#961988, David Kopriva).
Gregor Gassner and Andrés M. Rueda-Ramírez acknowledge funding through the Klaus-Tschira Stiftung via the project ``HiFiLab'' (00.014.2021) and through the German Federal Ministry for Education and Research (BMBF) project ``ICON-DG'' (01LK2315B) of the  ``WarmWorld Smarter'' program.
Gregor Gassner and Daniel Bach acknowledge funding from the German Science Foundation DFG through the research unit ``SNuBIC'' (DFG-FOR5409).

Andrés Rueda-Ramírez gratefully acknowledges funding from the Spanish Ministry of Science, Innovation, and Universities through the ``Beatriz Galindo'' grant (BG23-00062).
\bibliography{references}

\begin{thebibliography}{8}
\expandafter\ifx\csname natexlab\endcsname\relax\def\natexlab#1{#1}\fi
\providecommand{\url}[1]{\texttt{#1}}
\providecommand{\href}[2]{#2}
\providecommand{\path}[1]{#1}
\providecommand{\DOIprefix}{doi:}
\providecommand{\ArXivprefix}{arXiv:}
\providecommand{\URLprefix}{URL: }
\providecommand{\Pubmedprefix}{pmid:}
\providecommand{\doi}[1]{\href{http://dx.doi.org/#1}{\path{#1}}}
\providecommand{\Pubmed}[1]{\href{pmid:#1}{\path{#1}}}
\providecommand{\bibinfo}[2]{#2}
\ifx\xfnm\relax \def\xfnm[#1]{\unskip,\space#1}\fi
%Type = Article
\bibitem[{Arnold~DN(2006)}]{feec}
\bibinfo{author}{Arnold~DN, Falk~RS, W.R.}, \bibinfo{year}{2006}.
\newblock \bibinfo{title}{{F}inite element exterior calculus, homological techniques, and applications}.
\newblock \bibinfo{journal}{Acta Numerica} \bibinfo{volume}{15}, \bibinfo{pages}{1--155}.
\newblock \DOIprefix\doi{https://doi.org/10.1017/S0962492906210018}.
%Type = Inproceedings
\bibitem[{Gerritsma(2011)}]{gerritsma}
\bibinfo{author}{Gerritsma, M.}, \bibinfo{year}{2011}.
\newblock \bibinfo{title}{{E}dge {F}unctions for {S}pectral {E}lement {M}ethods}, in: \bibinfo{editor}{Hesthaven, J.S.}, \bibinfo{editor}{R{\o}nquist, E.M.} (Eds.), \bibinfo{booktitle}{{S}pectral and {H}igh {O}rder {M}ethods for {P}artial {D}ifferential {E}quations}, \bibinfo{publisher}{Springer Berlin Heidelberg}, \bibinfo{address}{Berlin, Heidelberg}. pp. \bibinfo{pages}{199--207}.
\newblock \DOIprefix\doi{https://doi.org/10.1007/978-3-642-15337-2_17}.
%Type = Article
\bibitem[{Jasper~Kreeft(2011)}]{gerritsma_mimetic_framework}
\bibinfo{author}{Jasper~Kreeft, Artur~Palha, M.G.}, \bibinfo{year}{2011}.
\newblock \bibinfo{title}{{M}imetic framework on curvilinear quadrilaterals of arbitrary order} \DOIprefix\doi{https://doi.org/10.48550/arXiv.1111.4304}.
%Type = Article
\bibitem[{Kopriva(2006)}]{kopriva}
\bibinfo{author}{Kopriva, D.}, \bibinfo{year}{2006}.
\newblock \bibinfo{title}{{M}etric {I}dentities and the {D}iscontinuous {S}pectral {E}lement {M}ethod on {C}urvilinear {M}eshes}.
\newblock \bibinfo{journal}{Journal of Scientific Computing} \bibinfo{volume}{26}, \bibinfo{pages}{301--327}.
\newblock \DOIprefix\doi{https://doi.org/10.1007/s10915-005-9070-8}.
%Type = Inproceedings
\bibitem[{Mark H.~Carpenter(1994)}]{carpenter}
\bibinfo{author}{Mark H.~Carpenter, C.A.K.}, \bibinfo{year}{1994}.
\newblock \bibinfo{title}{{F}ourth-{O}rder 2n-{S}torage {R}unge- {K}utta {S}chemes}.
%Type = Article
\bibitem[{Nédélec(1980)}]{nedelec}
\bibinfo{author}{Nédélec, J.C.}, \bibinfo{year}{1980}.
\newblock \bibinfo{title}{{M}ixed finite elements in $\mathbb{R}^3$}.
\newblock \bibinfo{journal}{Numerische Mathematik} .
%Type = Article
\bibitem[{Ranocha et~al.(2022)Ranocha, Schlottke-Lakemper, Winters, Faulhaber, Chan and Gassner}]{ranocha2022adaptive}
\bibinfo{author}{Ranocha, H.}, \bibinfo{author}{Schlottke-Lakemper, M.}, \bibinfo{author}{Winters, A.R.}, \bibinfo{author}{Faulhaber, E.}, \bibinfo{author}{Chan, J.}, \bibinfo{author}{Gassner, G.}, \bibinfo{year}{2022}.
\newblock \bibinfo{title}{Adaptive numerical simulations with {T}rixi.jl: {A} case study of {J}ulia for scientific computing}.
\newblock \bibinfo{journal}{Proceedings of the JuliaCon Conferences} \bibinfo{volume}{1}, \bibinfo{pages}{77}.
\newblock \DOIprefix\doi{https://doi.org/10.21105/jcon.00077}, \href{http://arxiv.org/abs/2108.06476}{\tt arXiv:2108.06476}.
%Type = Misc
\bibitem[{Schlottke-Lakemper et~al.(2021)Schlottke-Lakemper, Gassner, Ranocha, Winters and Chan}]{schlottkelakemper2020trixi}
\bibinfo{author}{Schlottke-Lakemper, M.}, \bibinfo{author}{Gassner, G.J.}, \bibinfo{author}{Ranocha, H.}, \bibinfo{author}{Winters, A.R.}, \bibinfo{author}{Chan, J.}, \bibinfo{year}{2021}.
\newblock \bibinfo{title}{{T}rixi.jl: {A}daptive high-order numerical simulations of hyperbolic {PDE}s in {J}ulia}.
\newblock \bibinfo{howpublished}{\url{https://github.com/trixi-framework/Trixi.jl}}.
\newblock \DOIprefix\doi{https://doi.org/10.5281/zenodo.3996439}.

\end{thebibliography}
\bibliographystyle{cas-model2-names}

\end{document}